\documentclass[11pt]{amsart}
\usepackage{amsmath,amssymb, xypic,verbatim,amscd}
\numberwithin{equation}{section}

\newcommand\bull{\sssize{\bullet}}
\newcommand\pto{\rho}
\newcommand\deltaplus{\Delta_+}
\newcommand\deltaminus{\Delta_-}
\newcommand\D{\mathcal{D}}
\newcommand\tv{\tilde{V}}
\newcommand\ttau{\tilde{\tau}}
\newcommand\ttheta{\tilde{\theta}}
\newcommand\tq{\tilde{Q}}

\newcommand\mvv{V}

\newcommand\Ri{R(U(r))_{k,n}}
\newcommand\QH{\operatorname{QH}^*}
\newcommand\mqq{Q}
\newcommand\mee{E}
\newcommand\mff{F}
\newcommand\Ml{\mathcal{L}}

\newcommand\mii{I}
\newcommand\mjj{J}
\newcommand\mll{L}
\newcommand\spane{\operatorname{Span}}
\newcommand\home{\operatorname{Hom}}
\newcommand\uhom{\underline{\operatorname{Hom}}}
\newcommand\spec{\operatorname{Spec}}

\newcommand\mf{\mathcal{F}}

\newcommand\pone{\Bbb{P}^1}

\newcommand\td{\tilde{d}}

\newcommand\SU{SU}

\newcommand\tensor{\otimes}
\newcommand\ml{\mathcal{L}}

\newcommand\ma{\mathcal{A}}
\newcommand\shom{\mathcal{H}{\operatorname{om}}}
\newcommand\codim{\operatorname{codim}}

\newcommand\GL{\operatorname{GL}}

\newcommand\im{\operatorname{im}}

\newcommand\rk{\operatorname{rk}}

\newcommand\Gr{\operatorname{Gr}}
\newcommand\SL{\operatorname{SL}}

\newcommand{\leto}[1]{\stackrel{#1}{\to}}

\newcommand\mfm{\operatorname{Gr}}

\newcommand\mpp{{\operatorname{ S}}}

\newtheorem{theorem}{Theorem}[section]

\newtheorem{proposition}[theorem]{Proposition}

\newtheorem{lemma}[theorem]{Lemma}

\newtheorem{defi}[theorem]{Definition}

\newtheorem{remark}[theorem]{{ Remark}}

\begin{document}
\title[Duality theorem]{The strange duality conjecture for generic curves}
\author{Prakash Belkale}\footnote{The author was partially
supported by NSF grant DMS-0300356.}
\address{Department of Mathematics\\ UNC-Chapel Hill\\ CB \#3250, Phillips Hall
\\ Chapel Hill,  NC 27599}
\email{belkale@email.unc.edu}
\begin{abstract}
Let $X$ be a smooth connected projective algebraic curve of genus
$g\geq 1$. The strange duality conjecture connects
 non-abelian theta functions of rank $r$ and level $k$, and those of rank $k$ and
level $r$ on $X$ (for $\SU(r)$ and $\operatorname{U}(k)$
respectively). In this paper we prove this conjecture for $X$
generic in the moduli space of curves of genus $g$.
\end{abstract}

 \maketitle
\section{Introduction}
Let $SU_X(r)$ be the moduli space of semi-stable vector bundles of
rank $r$ with trivial determinant over a connected smooth projective
algebraic curve $X$ of genus $g\geq 1$ over $\Bbb{C}$. Recall that a
vector bundle $E$ on $X$ is called semi-stable  if for any subbundle
$V$, $\deg(V)/\rk(V) \leq \deg(E)/\rk(E)$. Points of $\SU_X(r)$
correspond to isomorphism classes of semi-stable rank $r$ vector
bundles  with trivial determinant up to an equivalence relation.

For any line bundle $L$ of degree $g-1$ on $X$ define
$\Theta_L=\{E\in \SU(r),h^0(E\tensor L)\geq 1\}$. This turns out be
a non-zero Cartier divisor whose associated line bundle
$\ml=\mathcal{O}(\Theta_L)$ does not depend upon $L$. It is known
that $\ml$   generates the Picard group of $\SU_X(r)$ (for this and
the  precise definition of $\ml$ in terms of determinant of
cohomology  see ~\cite{DN}).

Let $U_X^*(k)$ be the moduli space of semi-stable rank $k$ and
degree $k(g-1)$ bundles on $X$. Recall that on $U_X^*(k)$ there is a
canonical non-zero theta (Cartier) divisor $\Theta_k$ whose
underlying set is $\{F\in U_X^*(k), h^0(X,F)\neq 0 \}.$ Put
$\mathcal{M}=\mathcal{O}(\Theta_k)$. Consider the natural map
$\tau_{k,r}:\SU_X(r)\times U^*_X(k)\to U_X^*(kr)$ given by tensor
product. From the theorem of the square, it follows
 that $\tau_{k,r}^*\mathcal{M}$ is isomorphic to $\ml^k\boxtimes \mathcal{M}^r$. The canonical element $\Theta_{kr}\in
H^0(U_X^*(kr),\mathcal{M})$ and the Kunneth theorem gives a map
 well defined up to scalars:
\begin{equation}\nonumber\tag{$\dagger$}
H^0(U_X^*(k),\mathcal{M}^r)^*\to H^0(\SU_X(r),\mathcal{L}^k).
\end{equation}
The strange duality conjecture asserts that ($\dagger$) is an
isomorphism. It is known that $h^0(U_X^*(k),\mathcal{M}^r)$ equals
$h^0(\SU_X(r),\mathcal{L}^k)$ (see for example ~\cite{beau}, Section
8). The strange duality conjecture is known to hold when
\begin{itemize}
\item $k=1$, and arbitrary
$r$ ~\cite{bnr}.
\item $r=2$, $k=2$ and $C$ has no vanishing thetanull ~\cite{BBB}.
\item $r=2$, $k=4$ and $C$ has no vanishing thetanull ~\cite{vGP}.
\item $r=2$, $k$ even and $k\geq 2g-4$, and generic $C$ ~\cite{L}.
\end{itemize}

 An element $F$ of $U_X^*(k)$ produces an element
$\Theta_F$ of $H^0(SU_X(r), \ml^k)$  well defined up to scalars. The
zero locus of $\Theta_F$  is the set of all $E\in\SU_X(r)$ such that
$h^0(X,E\tensor F)\neq 0$ (the degree of $F$ is such that
$\chi(X,E\tensor F))=0$). It is easy to see that ($\dagger$) is an
isomorphism if and only if $\Theta_F$ for $F\in U_X^*(k)$ span
$H^0(\SU_X(r), \ml^k)$.

Let $M_g$ denote the moduli-space of connected smooth projective
algebraic curves of genus $g$. In this paper, we prove the strange
duality conjecture for generic curves:
\begin{theorem}\label{al}
For generic $X\in M_g$, the sections $\Theta_F$ for $F\in U_X^*(k)$
span $H^0(SU_X(r), \Ml^k)$.
\end{theorem}
 For the history of this problem as well as recent developments we
refer the reader to Beauville ~\cite{beau}, Donagi-Tu
~\cite{donagitu}, Polishchuk ~\cite{polish} and  Popa ~\cite{popa}.

 \subsection{The main idea} The starting point for this
paper is the classical relation between the cohomology of
Grassmannians $\Gr(r,n)$ and invariant theory of the special linear
group $\SL(r)$ (or equivalently $\SU(r)$). In ~\cite{invariant},
this relation was further strengthened by demonstrating how triple
intersections of Schubert varieties geometrically produce a basis
for the invariants for the associated $\SL(r)$ tensor product
multiplicity problem. The next step is Witten's ~\cite{witten}
relation between the (small) quantum cohomology of Grassmannians and
structure coefficients in the Verlinde algebra for $\SU(r)$ (these
are dimensions of spaces of sections of theta bundles on moduli
spaces of parabolic bundles on $\pone$). This relation can be
geometrized in a similar way. Theorem ~\ref{al} of this paper is a
higher genus generalization of this relation with
$H^0(SU_X(r),\ml^k)$ viewed as a representation theoretic object.
The associated enumerative problem needs to be {\em invented}. The
linearly independent sections coming from the inherent
transversality in the enumerative problem will be shown to have the
form $\Theta_F$ for $F\in U_X^*(k)$.

To help us invent the enumerative problem that should correspond
 to $H^0(SU_X(r), \Ml^k)$, we calculate the dimension $M(r,k,g)$ of the latter (it is known that the rank of $H^0(SU_X(r), \Ml^k)$
  does not vary with $X\in M_g$).
  Using a factorization formula, the dimension can be related to the
dimensions of conformal blocks for $\pone$. This reduction uses the
Verlinde formula (Beauville-Laszlo ~\cite{v2}, Faltings ~\cite{v1}
and Kumar-Narasimhan-Ramanathan ~\cite{v3}), and the factorization
formula of Tsuchiya-Ueno-Yamada ~\cite{TUY}. Now, by a theorem of
Witten (see ~\cite{witten}, a mathematical proof was given by
Agnihotri ~\cite{agni})
 there is a relation between conformal blocks for $\pone$ and  the  (small) quantum cohomology of
 Grassmannians. Putting all these together one finds a formula (see Section ~\ref{examples} for some examples):
\begin{equation}\nonumber\tag{$\ddagger$}
M(r,k,g)=\sum
 '\langle\omega_{I^1},\dots,\omega_{I^g},\omega_{(I^1)'},\dots,\omega_{(I^g)'}\rangle_{0,-k(g-1)}
\end{equation}
where the sum is over all sequences of subsets $(I^1,\dots,I^g)$ of
$[r+k]=\{1,2,\dots,r+k\}$ with $r$ elements, {\em each of which
contain $1$}, and  where the Gromov-Witten invariants are
``twisted'' (see Section ~\ref{temina} for the definition, and
Section ~\ref{fultone} for the definition of the classes
$\omega_I$).

The simultaneous appearance  of Schubert cohomology classes
$\omega_I$ and their duals $\omega_{I'}$ in Equation ($\ddagger$)
leads one to suspect the role of the diagonal in a product of
Grassmannians $\Gr(r,n)\times \Gr(r,n)$. However, the restricted
nature of the sum suggests a piece in a partial degeneration of the
diagonal ($\deltaplus $ in Section ~\ref{godfather}).

Using the insight given to us by Equation ($\ddagger$), we introduce
an enumerative problem that corresponds to $M(r,k,g)$ (see Section
~\ref{before}).

 After this work was completed, I received a preprint
from Takashe Abe ~\cite{abe}  in which he proves the  strange
duality conjecture for generic curves, for $r=2$, arbitrary $k$. In
these cases
 he  proves the stronger version of strange duality (generalized degrees)
 formulated in Donagi-Tu ~\cite{donagitu}.

I thank A. Boysal, P. Brosnan, M. V. Nori and M. Popa for  useful
discussions. M. Popa helped me correct an initial naive idea.
\section{Determinant of the cohomology}\label{falti}
Let $\mathcal{L}$ be a line bundle on a (possibly non complete)
variety $Z$ and $s_1,\dots,s_m$ global sections of $\ml$. One way of
proving that $s_1,\dots,s_m$ are linearly independent is to find
points $z_1,\dots,z_m$ in $Z$ so that the determinant of the matrix
$(s_i(z_j))$, $i=1,\dots,m$ and $j=1,\dots,m$ is non-zero. To make
sense of this determinant note that $s_i(z_j)\in \ml_{z_j}$ and
hence the determinant is an  element of $\tensor_{j=1}^m \ml_{z_j}$.
We view $z_j$ as witnesses to the linear independence of the $s_i$.

In this section we formulate  the notion of an ``$(m,r,k)$-frame''.
Intuitively an $(m,r,k)$-frame $(E_1,\dots,E_{m};F_1,\dots,F_{m})$
on a possibly singular irreducible curve $C$ of arithmetic genus $g$
is
 a tuple $(E_1,\dots,E_{m};F_1,\dots,F_{m})$ where the $E_j$ are vector bundles on $C$ of
rank $r$ degree $0$
 and with isomorphic determinantal line bundles; $F_i$ are vector bundles on $C$   of rank $k$ and degree
$k(g-1)$, so that  assuming $\det(E_i)\leto{\sim} \mathcal{O}$
(otherwise we will have to twist by a line bundle of degree 0),
$E_j\in \SU_X(r)$'s are witnesses to the linear independence of the
sections $\Theta_{F_i}$. We will need to clarify the meaning of
theta functions if $C$ is singular (we do this without compactifying
the moduli of bundles on $C$).

The reader may now wish to turn to Section ~\ref{nodal} on a first
reading.
\subsection{Determinant of cohomology and basic operations}\label{falt}
 Let $\pi:X\to S$ be a relative curve (that is, it is flat, the
fibers are proper of dimension $\leq 1$ and
$\pi_*(\mathcal{O}_X)=\mathcal{O}_S$). If $\mff$ is a coherent sheaf
on $X$ which is flat over $S$, we form the determinant of its
cohomology $\mathcal{D}(\mff)$ which  is a line bundle on $S$
(following ~\cite{EGA3}, Section 7). Formally, the fiber of
$\D(\mff)$ at a point $s\in S$ is the one dimensional vector space
$$\det H^0(X_s,\mff_s)^*\tensor
\det H^1(X_s,\mff_s).$$ More precisely, the push forward in the
 bounded derived category $R\pi_*(\mff)$ is represented locally by a complex
$\mf_0\to \mf_1$ of vector bundles. $\D(\mff)$ is defined to be
$\det \mf^*_0\tensor \det \mf_1$ (where $\det \mf$ denotes the top
exterior power of the vector bundle $\mf$). This local definition
globalizes.

 Suppose that  $\forall s\in S$, $\chi(\mff_s)=0$. Then, $R\pi_*(\mff)$ is locally represented
 by a complex $\mf_0\leto{\psi} \mf_1$ with $\rk(\mf_0)=\rk(\mf_1)$.
Taking the top exterior product of $\psi$, we find a canonical
section $\sigma(\mff)\in \D(\mff)$. It is easy to see that for $s\in
S$, $\sigma(F_s)=0$ if and only if $h^0(X_s,F_s)\neq 0.$

The determinant of cohomology and its canonical section satisfy some
compatibility properties:
\begin{enumerate}
\item If $\alpha: 0\to \mff_1\to\mff_2\to\mff_3\to 0$ is an exact
sequence of $S$-flat coherent sheaves on $X$, then there is an
induced isomorphism
$$\D(\alpha):\D(\mff_2)\leto{\sim}\D(\mff_1)\otimes\D(\mff_3).$$
\item If in (1), the relative Euler characteristics of
$\mff_1$, $\mff_2$ and $\mff_3$ are each zero, then
$$\D(\alpha)(\sigma(\mff_2))=\sigma(\mff_1)\otimes\sigma(\mff_3).$$
\end{enumerate}
In ~\cite{faltings}, G. Faltings  makes the following definition
\begin{defi}
Let $I$ and $J$ be $S$-flat coherent sheaves on $X$. We say $I$ and
$J$ coincide generically in K-theory if the following holds: For all
$s\in S$ such that $\operatorname{depth}(\mathcal{O}_{S,s})=0$, the
difference of $I$ and $J$ in the K-theory of coherent
$\mathcal{O}_{S,s}$ flat sheaves on $X\times_S
\operatorname{Spec}(\mathcal{O}_{S,s})$ is represented by a finite
complex of sheaves, such that the support of  its cohomology  is
finite over $\spec{\mathcal{O}}_{S,s}$.
\end{defi}
If, for example $X$ and $S$ are reduced and irreducible, and $I$ and
$J$ are $S$-flat and have the same generic rank (that is, rank at
the generic point of $X$), then they coincide generically in
K-theory. This case will be sufficient for this paper.

Faltings then  proves (Theorem I.1 of ~\cite{faltings}) that if
${I}$ and ${J}$ are two $S$-flat coherent sheaves on $X$ of the same
relative Euler characteristic
 and which coincide generically in $K$-theory, and $E$, $E'$ vector
bundles on $X$ of the same rank with a given isomorphism
$\psi:\det(\mee)\to\det(E')$, there is a natural isomorphism
$$\Phi_{\mee,E'}(\psi):\D(E\tensor I)\tensor\D(E'\tensor J)\leto{\sim}\D(E'\tensor I)\tensor\D(E\tensor J)$$
with various functorial and compatibility properties (compatible
exact sequences in $(\mee,E')$). {\em We will always assume that $I$
and $J$ are in addition also of the same Euler characteristic}. This
ensures that $\Phi_{E,E'}(\psi)$ does not depend upon the choice of
the isomorphism $\psi:\det(E)\leto{\sim}\det(E')$. We will therefore
drop the dependence on $\psi$ in the maps $\Phi_{E,E'}$.

 In the remainder
of this section we will recall the construction of $\Phi_{E,F}$ and
 note its basic properties.

Consider the special case when exist line bundles
$\mll_1,\dots,\mll_k$ and filtrations
$$\mee=\mee_1\supset\mee_2\supset\dots\supset\mee_r\supset 0=\mee_{r+1},$$
 $$E'=E'_1\supset E'_2\supset\dots\supset E'_r\supset 0=E'_{r+1}$$
and (given) isomorphisms $\mee_j/\mee_{j+1}\to L_j$ and
$E'_j/{E'_{j+1}}\to L_j$ for $j=1,\dots, r$ . Such a situation can
be arranged after a suitable flat base change of $S$  (for details
see ~\cite{faltings}). If $S=\operatorname{Spec}(\Bbb{C})$, no base
changes are needed!. In such a set up there are induced isomorphisms
$$A:\D(E\tensor I)\to \prod_{j=1}^r\D(\mll_j\tensor I),\ B:\D(E'\tensor J)\to \prod_{j=1}^r\D(\mll_j\tensor J).$$
$$C:\D(E'\tensor I)\to \prod_{j=1}^r\D(\mll_j\tensor I),\ D:\D(E\tensor J)\to \prod_{j=1}^r\D(\mll_j\tensor J)$$
The isomorphism $\Phi_{\mee,E'}$ in ``this good coordinate system''
is just the identity map. That is the diagram below commutes (where
the bottom arrow is identity):
$$\xymatrix{\D(\mee\tensor\mii)\tensor\D(E'\tensor J) \ar[d]^{A\tensor B}\ar[r]^{\Phi_{\mee,E'}} & \D(E'\tensor \mii)\tensor\D(E\tensor\mjj)\ar[d]^{C\tensor D}\\
\prod_{j=1}^r\D(\mll_j\tensor\mii)\tensor\prod_{j=1}^r\D(\mll_j\tensor
J)\ar[r]^{=} & \prod_{j=1}^r\D(\mll_j\tensor
I)\tensor\prod_{j=1}^r\D(\mll_j\tensor J)
\\}$$

Faltings shows that $\Phi_{E,E'}$ is independent of the good
coordinates (by this we mean the choice of filtrations of $E$ and
$E'$). Together with ~\cite{emperor}, Proposition 1 (ii), it is easy
to see that $\Phi_{E,E'}$ is compatible with exact sequences in
$(I,J)$. That is, if
$$0\to I_1\to I\to I_2\to 0\text{, } \ 0\to J_1\to J\to J_2\to 0$$ are exact sequences of $S$-flat
 coherent sheaves such that $I_a,J_a$ have the same relative Euler
characteristics for $a=1,2$ and coincide generically in $K$-theory,
then the following diagram commutes ($E$ and $E'$ are as above)
$$\xymatrix@R=0.05cm{\D(E\tensor I)\tensor \D(E'\tensor
J)\ar[r]^{\Phi_{E,E'}}\ar[dd] & \D( E'\tensor I)\tensor \D(E\tensor
J)\ar[dd] \\
\\
D(E\tensor I_1)\tensor \D(E'\tensor
J_1)\ar[r]\tensor\ar[r]^{\Phi_{E,E'}} & D(E'\tensor I_1)\tensor
\D(E\tensor J_1)\tensor\\
D(E\tensor I_2)\tensor \D(E'\tensor J_2) &  D(E'\tensor I_2)\tensor
\D(E\tensor J_2)\\
}$$ where the bottom $\Phi_{E,E'}$ is the tensor product of
$\Phi_{E,E'}$ acting on the two levels. Faltings does not list this
property, and we learned it from the proof of Lemma 1 in
~\cite{este}.

The second observation, which is obvious in good coordinates, is
that if $\mii=\mjj$, then the isomorphism $\Phi_{E,E'}$ is just the
permutation of the two factors (if $L$ is a one dimensional vector
space $\ell_1,\ell_2\in L$, then $\ell_1\tensor
\ell_2=\ell_2\tensor\ell_1\in L\tensor L$).

The third observation is that if $(E_1,\dots,E_m;F_1,\dots,F_m)$ are
vector bundles on $X$ such that  $E_1,\dots, E_m$ are vector bundles
with isomorphic determinants and  each of rank $r$, $F_1,\dots,F_m$
are vector bundles on $X$  of the same relative degree and rank; and
$\pi\in S_n$ (the symmetric group on $\{1,\dots,n\}$) then there are
 ``natural''isomorphisms, compatible with compositions:
$$\D(E_1\tensor F_1)\tensor\D(E_2\tensor
F_2)\tensor\dots\tensor\D(E_m\tensor F_m)\leto{A_{\pi}}
\D(E_{\pi(1)}\tensor F_1)\tensor\D(E_{\pi(2)}\tensor
F_{2})\tensor\dots\tensor\D(E_{\pi(m)}\tensor F_{m}).$$ To see this,
there are such maps for transpositions $\pi$ viz. $\Phi_{E_i,E_j}$.
Writing any permutation as a composition of transpositions, we will
obtain maps for all $\pi$. The independence from the choice of the
representation of $\pi$ as a composition of permutations (as well as
compatibility properties) follows from: ``in good coordinates for
$E_i$'' all these maps are identity.
\subsection{ Frames of vector bundles on curves}\label{lecteur}
 Let $C$ be a reduced irreducible projective algebraic curve  of arithmetic genus
$g$ with only ordinary double points for singularities. Consider a
tuple $(E_1,\dots,E_m;F_1,\dots,F_m)$ where
\begin{enumerate}
\item[(P1)] $E_1,\dots, E_m$ are
vector bundles on $C$ each of degree $0$, rank $r$ and such that for
$i,j\in\{1,\dots,m\}$, $\det(E_i)$ is isomorphic to $\det(E_j)$.
\item[(P2)] $F_1,\dots,F_m$ are vector bundles on $C$, each of rank $k$ and degree $k(g-1)$.
\end{enumerate}
Notice that $\chi(C,E_i\tensor F_j)=0$ and hence we have canonical
elements $\sigma(E_i\tensor F_j)\in \D(E_i\tensor F_j)$ (a one
dimensional vector space).

By the third observation in Section ~\ref{falt},  the one
dimensional vector spaces for various $\pi\in S_n$ (the symmetric
group)
$$
\D(E_{\pi(1)}\tensor F_1)\tensor\D(E_{\pi(2)}\tensor
F_{2})\tensor\dots\tensor\D(E_{\pi(m)}\tensor F_{m})$$
 are all
canonically identified. Consider
 the $m\times m$ matrix $(\sigma(E_i\tensor F_j))$,
$i=1,\dots,m$, $j=1,\dots,m$. Because of the remark above it makes
sense to speak about the non-zeroness of the determinant of this
matrix.
\begin{defi}\label{framedefi}
 A tuple $(E_1,\dots,E_m;F_1,\dots,F_m)$ satisfying properties $(P1)$ and $(P2)$ is said to be an $(m,r,k)$-frame on $C$ if the determinant
of the $m\times m$ matrix $(\sigma(E_i\tensor F_j))$, $i=1,\dots,m$,
$j=1,\dots,m$ is non-zero.
\end{defi}
We will give a geometric interpretation of the definition of a
frame. Let $(E_1,\dots,E_m;F_1,\dots,F_m)$ be a tuple as above. Fix
a vector bundle $E_0$ of the same rank as $E_i$'s whose determinant
is isomorphic to that of $E_i$ (for example one of the $E_i$). Let
$F_0$ be a vector bundle of the same degree and rank as the $F_j$'s.
Choose nonzero elements $s_j\in \D(E_0\otimes F_j)$ and
$t\in\D(E_0\tensor F_0)$. These choices do not involve
$E_1,\dots,E_m$ . We have canonical isomorphisms
$$\Phi_{E,E_0}:\D(E\tensor F_i)\tensor\D(E_0\tensor F_0)\leto{\sim}\D(E_0\tensor F_i)\tensor\D(E\tensor F_0)$$
The choice of $t,s_0,\dots,s_m$ gives a morphism
$\lambda_j:\D(E\tensor F_j)\to\D(E\tensor F_0)$.

Therefore one obtains $m$ tuples of complex numbers well defined up
to non zero scalars  $\Theta(E)=[\lambda_1(\sigma(E\tensor
F_1)),\dots,\lambda_m(\sigma(E\tensor F_m))]$. But $\Theta(E)$
 (up to non-zero scalars) does depend upon the choice of
 $(E_0,F_0,t,s_1,\dots,s_m)$. For example, if we change $s_i$, then
 the entry $i$ of $\Theta(E)$ gets scaled.

It is easy to verify that $(E_1,\dots,E_m; F_1,\dots,F_m)$ is a
$(m,r,k)$-frame if and only $\Theta(E_1),\dots,\Theta(E_m)$ are
linearly independent.

 The main properties of $(m,r,k)$-frames are:

\begin{proposition}\label{map}
Let  $X\in M_g$  and $(E_1,\dots,E_m;F_1,\dots,F_m)$ an $(m,r,k)$
frame on $X$ with $\det(E_i)\leto{\sim}\mathcal{O}_X$. Then
$\Theta_{F_j}\in H^0(\SU_X(r),\Ml^k),$ $j=1,\dots,m$ are linearly
independent. If $m\geq\rk H^0(\SU_X(r),\Ml^k)$, these sections give
a basis of $H^0(\SU_X(r),\Ml^k)$ (and hence $m=\rk
H^0(\SU_X(r),\Ml^k)$.)
\end{proposition}
\begin{proof}
From the determinantal condition in the definition of the frame, we
see that for each $i$ there exists a $j$ (similarly for each $j$
there is an $i$) so that $\sigma(E_i\tensor F_j)\neq 0$. Therefore
$E_i$ and $F_j$ are semi-stable vector bundles on $X$ for all $i$
and $j$ (Use Lemma ~\ref{nori}). For the linear independence we can
replace $\SU_X(r)$ by a space $Z$ that carries a universal family of
semistable vector bundles of trivial determinant and rank $r$ on
$X$. The same determinantal condition now implies the linear
independence of $\Theta_{F_j}$.
\end{proof}
\begin{lemma}\label{don}
Suppose that a $(m,r,k)$- frame $(E_1,\dots,E_m;F_1,\dots,F_m)$
 exists on an reduced irreducible projective curve $C$ with at most double point singularities . Then there
 exists such an $(m,r,k)$-frame on $C$ with $\det(E_i)$ isomorphic to $\mathcal{O}_C$.
\end{lemma}
\begin{proof}
The group of line bundles of  degree $0$  on $C$  is a divisible
group. Suppose that $L^r\leto{\sim}\det(E_i)$ (recall that we are
assuming that $E_i$ is of degree $0$). Then we consider the modified
frame
$$(E_1\tensor L^{-1},\dots, E_m\tensor
 L^{-1};F_1\tensor L,\dots, F_m\tensor L)$$
\end{proof}
\begin{lemma}\label{moment}Consider a projective flat family  $X\to S$ of
curves. Let $s_0\in S$ be such that $X_{s_0}$ is a reduced
irreducible curve with at most ordinary double point singularities .
Suppose that a $(m,r,k)$- frame
 exists on $C_{s_0}$.  Then there is an  open subset $U\subseteq S$ containing $s_0$ such that $X_s$ has an
 $(m,r,k)$-frame  for each $s\in U$.
\end{lemma}
\begin{proof} Let $(E_1,\dots,E_m;F_1,\dots,F_m)$ be an
$(m,r,k)$-frame on $X_{s_0}$. Using Lemma ~\ref{don} assume that
$\det(E_i)\leto{\sim}\mathcal{O}_{X_{s_0}}$ for $i=1,\dots,m$.

Using Lemma ~\ref{deform}, find an  \'{e}tale   base change $S'\to
S$ so that the image contains $s_0$ and $(E_1\dots,E_m;
F_1,\dots,F_m)$ lift to vector bundles on $C_{S'}$ satisfying
properties (a) and (b) of the frame. The non vanishing of the
determinant is an open condition and this concludes the proof.
\end{proof}

\section{Line bundles and vector bundles on   rational nodal
curves}\label{nodal}
 Let us fix a rational nodal curve $C$ of arithmetic genus $g$ together with a
normalization $\pone\leto{f} C$. Let $r_1,\dots,r_g$ be the nodes of
$C$ and let $f^{-1}(r_j)=\{p_j,q_j\}$ for $j=1,\dots,g$. This
notation will be fixed throughout the paper.

Let $\mll$
 be a line bundle of degree $0$ on $C$ (that is, the degree of $f^*L$ is zero). It is clear that $f^*\mll$ is trivial, that is there is a
  isomorphism unique up to scalars $A:\mathcal{O}_{\pone}\to f^*\mll$.
  The canonical map $f^*(\mll)_{p_j}\to f^*(\mll)_{q_j}$ gives us a
  well defined scalar $c_j$ which makes the following diagram commute:
  $$\xymatrix{
  \Bbb{C}\ar[d]^{A_{p_j}}\ar[r]^{\cdot c_j} & \Bbb{C}\ar[d]^{A_{q_j}} \\
  f^*\mll_{p_j}\ar[r]^{=}  &  f^*\mll_{q_j}}$$
We therefore obtain a well defined morphism of groups
$$\{\text{Line bundles of degree 0 on C}\}\  \to \  (\Bbb{C}^*)^g,\ \mll\leadsto (c_1,\dots,c_g).$$
It is easily checked that this is a group isomorphism.
\subsection{$S$-bundles on rational nodal curves}
\begin{defi}\label{hai}
An $S$-bundle (S for strange) on $C$ is a pair $(\mvv,\tau)$ where
$V$ is a vector bundle on $\pone$ together with a $g$-tuple
$(\tau_1,\dots,\tau_g)$ where $\tau_{j}:\mvv_{p_j}\to \mvv_{q_j}$
for $j=1,\dots,g$ are morphisms of vector spaces.
\end{defi}
Consider the {\em surjective} morphism $f_{*}\mvv\to
\oplus_{j}\mvv_{q_j}|r_j$ ($\mvv_{q_j}|r_j$ is the skyscraper sheaf
supported at $r_j$ with fiber $\mvv_{q_j}$) corresponding to the map
$\mvv_{p_j}\oplus\mvv_{q_j}\to \mvv_{q_j}$ given by $-\tau_{j}$ on
the first factor and identity on the second. Let $\tilde{V}$ be the
kernel. It is a coherent sheaf on $C$, and will be called {\em the
coherent sheaf underlying} $(V,\tau)$. It is easy to see that if the
$\tau_{j}$ are isomorphisms for $j=1,\dots,g$ then $\tilde{V}$ is a
vector bundle on $C$ (in this case $\tilde{V}$ is the vector bundle
on $C$ obtained by gluing $V_{p_j}$ and $V_{q_j}$ via the map
$\tau_j$ for $j=1,\dots,g$.)

Let $(\mvv,\tau)$ and $(\mqq,\theta)$ be $S$-bundles as above. A
morphism  $(\mvv,\tau)\to (\mqq,\theta)$ is a morphism
$\Gamma:\mvv\to \mqq$ on $\pone$ such that for $j=1,\dots,g$ the
diagram

$$\xymatrix{\mvv_{p_j}\ar[r]^{\tau_j}\ar[d]^{\Gamma_{p_j}}  &                \mvv_{q_j}\ar[d]^{\Gamma_{q_j}}\\
            \mqq_{p_j}\ar[r]^{\theta_j}            &                    \mqq_{q_j} }$$
commutes. The (vector) space of such homomorphisms is denoted by
$\uhom((\mvv,\tau),(\mqq,\theta))$.

Let $(V,\tau^{(1)})$ and $(V,\tau^{(2)})$ be $S$-bundles on $C$ with
$V=\mathcal{O}_{\pone}^{\oplus r}$ such that $\tau^{(i)}_j$ are
isomorphisms for $i=1,2$, $j=1,\dots,g$. These give rise to coherent
sheaves  $\tv_1$ and $\tv_2$ on $C$ which are locally  free sheaves.
We view $\tau^{(1)}_{j}$ and $\tau^{(2)}_{j}$ as $r\times r$
matrices using the identification
$V\leto{\sim}\mathcal{O}_{\pone}^{\oplus r}$. The following lemma is
immediate:
\begin{lemma}\label{looks}
As line bundles on $C$, $\det(\tv_1)$ and $\det(\tv_2)$ are
isomorphic if and only if
$\det(\tau^{(1)}_{j})=\det(\tau^{(2)}_{j})\in \Bbb{C}^*$ for
$j=1,\dots,g$.
\end{lemma}
\section{Outline of the argument}\label{before}
 Lemmas ~\ref{moment}, ~\ref{don} and Proposition
~\ref{map} give a strategy for proving the strange duality
conjecture generic $X\in M_g$.

{\bf Initial Strategy:} Construct an $(m,r,k)$-frame on the rational
nodal curve $C$ of arithmetic genus $g$ with $m\geq M(r,k,g)$
($M(r,k,g)$ was defined in the introduction).

 Here we use the fact that there exists  a family of
projective curves degenerating into $C$, such that the general
member of this family is a curve in $M_g$.

A $(m,r,k)$-frame on $C$ will be constructed by working on the
normalization $\pone$ of $C$. We will first recall some relevant
definitions and properties of ``evenly split'' bundles on $\pone$.
\subsection{Evenly split bundles on $\pone$}
A vector bundle $T$ on $\pone$ is said to be evenly split (ES) if
$W=\oplus_{i=1}^n\mathcal{O}_{\Bbb{P}^1}(a_i)$ with $\mid a_i-a_j
\mid \leq 1$ for $0<i<j\leq n$.

Let $D$, and $n$ be integers with $n>0$. It is easy to show there is
a unique, up to  isomorphism, ES -bundle of degree $-D$ and rank $n$
on $\pone$ (the negative sign is introduced for conformity with
notation used in quantum cohomology.)

 Let $T$ be a bundle on $\pone$. Define {$\mfm(d,r,T)$} to be the
moduli space of subbundles of $T$ which are of degree $-d$ and rank
$r$. This can be obtained as an open subset of the quot scheme of
quotients of $T$ of degree $d-D$ and rank $n-r$. In the notation of
~\cite{potier}, we are looking at the open subset of
$\operatorname{Hilb}^{n-r,d-D}(T)$ formed by points where the
quotient is locally free. The proof of the following standard result
may be found in ~\cite{qh} (Proposition 2.2).
\begin{proposition}\label{bigone}
Let $T$ be an ES bundle of degree $-D$ and rank $r$. Then,
$\mfm(d,r,T)$ is smooth and irreducible of dimension $r(n-r)+dn-Dr$.
Moreover, the subset of $\mfm(d,r,T)$ formed by ES-subbundles
$V\subseteq T$ such that $T/V$ is also ES, is open and dense in
$\mfm(d,r,T)$.
\end{proposition}
\subsection{The enumerative problem and the resulting frame on
$C$}\label{sharonstone}
 Let $T$ be a ES-vector bundle of rank
$n=r+k$ and degree $k(g-1)$ on $\pone$. Fix  a generic tuple
$$(\gamma_1,\dots,\gamma_g)\in\prod_{j=1}^g\home_{n-1}(T_{p_j},T_{q_j})$$
where $\home_{n-1}(T_{p_j},T_{q_j})$ denotes the set of maps of
vector spaces  $T_{p_j}\to T_{q_j}$ of rank $n-1$ (therefore the
kernel of $\gamma_j$ is one dimensional.)

The enumerative problem is the following: Count the number of
``singular S-subbundles'' $(V,\tau)$ of $(T,\gamma)$ so that
$\deg(V)=0$ and $\rk(V)=r$. More precisely, we want to count
subbundles $V$ of $T$ of degree $0$ and rank $r$ so that for
$j=1,\dots,g$,
\begin{itemize}
\item $\gamma_{j}(V_{p_j})\subset V_{q_j}$.
\item The induced map $\tau_j:V_{p_j}\to V_{q_j}$ is singular.
\end{itemize}
In Proposition ~\ref{shemoved}, we will show that there are only
finitely many such S-subbundles, and that scheme theoretically it is
reduced. The number $m$ of these will be shown to be $\geq M(r,k,g)$
(Proposition ~\ref{hotri}).

Denote them by $V^{(a)}$ for $a=1,\dots,m$. Let
$\mqq^{(a)}=T/V^{(a)}$ for $a=1,\dots,m$. Both $V^{(a)}$ and
$\mqq^{(a)}$ have natural S-bundle structures. We therefore obtain
S-bundles $(V^{(a)},\tau^{(a)})$ and $(Q^{(a)},\theta^{(a)})$. We
will show that the gluing maps $\theta^{(a)}_j$ for $\mqq^{(a)}$ are
nonsingular (Lemma ~\ref{versace}, (ii)). Hence the coherent sheaf
underlying $(\mqq^{(a)},\theta^{(a)})$ is a vector bundle
$\tilde{Q}^{(a)}$.

We perturb the maps $\tau^{(a)}$ (and leave maps $\theta^{(a)}$
unchanged) so that they become non-singular and such that the
corresponding vector bundles (on $C$) $\tilde{V}^{(a)}$ have
isomorphic determinant line bundles for $a=1,\dots,m$. This
perturbation is possible because we are starting from singular maps
$\tau^{(a)}_j$. Finally we will show, using the transversality in
the enumerative problem  that the perturbed tuple
$$((\tv^{(1)})^*,\dots,(\tv^{(m)})^*;\tq^{(1)},\dots,\tq^{(m)})$$
 is an $(m,r,k)$-frame on the  rational nodal curve $C$ (because in the limit an appropriate
 matrix is diagonal with non-zero entries on the diagonal). This will conclude the
 proof.
\section{Reformulation in terms of the diagonal}\label{godfather}
The enumerative problem given in Section ~\ref{sharonstone} is a
counting problem of number of ``singular'' subbundles of a bundle on
$C$ with singular gluing maps. However recall that if $\tau_j$ were
non singular then the set of subbundles of the induced bundle
$\tilde{T}$ (on $C$) of degree $0$ and rank $r$ can be described as
follows: Consider the natural map
$$\pi:\mfm(0,r,T)\to\prod_{j=1}^g\bigl( \Gr(r,T_{p_j}) \times \Gr(r,T_{q_j})\bigr)$$
Since there is a given isomorphism $\gamma_j:T_{p_j}\to T_{q_j}$ we
can define a ``diagonal''
$$\Delta_j= \{(A,B)\in \Gr(r,T_{p_j})\times \Gr(r,T_{q_j})\mid
\gamma_j(A)=B).$$ The set of subbundles of $\tilde{T}$ of degree $0$
and rank $r$ is just $\pi^{-1}\prod_{j=1}^g \Delta_j$.

In Section ~\ref{beethoven}, we show that the enumerative problem in
Section ~\ref{before}  has a similar description, but we need to
replace $\Delta_j$ by a singular ``diagonal''. To obtain
transversality in the enumerative problem and make use of it, we
will need to study tangent spaces as well.
\subsection{Degeneration of the diagonal}\label{degen}
Let $W$ be a vector spaces of rank $n$ and $0<r<n$ an integer. Let
$\Gr(r,W)$. We will consider partial degenerations of the diagonal
$\Delta\ =\{(A,B)\in \Gr(r,W)\times \Gr(r,W)\mid A \ = \  B\}$.

If $\Phi\in\operatorname{End}(W)$, then define
$$\Delta_{\Phi}=\{(A,B)\in\Gr(r,W)\times\Gr(r,W)\mid \Phi(A) \subseteq B\}.$$

We therefore obtain a subscheme
$$\tilde{\Delta}\subseteq \Gr(r,W)\times
\Gr(r,W)\times\operatorname{End}(W)$$  such that the fiber of the
projection  $T:\tilde{\Delta}\to \operatorname{End}(W)$ over $\Phi$
is $\Delta_{\Phi}$. Clearly $T$ is flat over
$\operatorname{Aut}(W)$. The map $\tilde{\Delta}\to
\Gr(r,W)\times\Gr(r,W)$  given by $(A,B,\Phi)\mapsto (A,B)$ is a
Zariski locally trivial fiber bundle with smooth fibers. Therefore
$\tilde{\Delta}$ is smooth.

We claim that $T$ is flat over endomorphisms of rank $\geq n-1$. It
suffices to show that if $\Phi$ is singular with kernel $L$ of rank
$1$, then $\Delta_{\Phi}$ is equidimensional of dimension $r(n-r)$.
If $(A,B)\in \Delta_{\Phi}$ then either $L\subset A$ or
$B\subset\im(\Phi)$. These two (irreducible) components  each have
the correct dimension (by a small calculation, for example the
second dimension is $\dim(\Gr(r,n-1)) +\dim(\Gr(r,r+1))$).

Let $L=\Bbb{C}\ell\subseteq W$ and $K\subseteq W$ be subspaces of
ranks $1$ and $n-1$ respectively so that the natural map $L\oplus K
\to W$ is an isomorphism. Let $\beta:W\to W$ be the corresponding
projection to $K$. $\Delta_{\beta}$ is a union
$\deltaplus\cup\deltaminus$ where
$$\deltaplus =\{(A,B)\in\Gr(r,W)\times\Gr(r,W)\mid \beta(A)\subseteq B, L\subseteq A\},$$
$$\deltaminus=\{(A,B)\in\Gr(r,W)\times\Gr(r,W)\mid \beta(A)\subseteq B, B\subseteq K\}.$$

We will need to be a bit more precise about the degeneration of
$\Delta$ to $\Delta_{\beta}$, and need a  more general degeneration
later on in Section ~\ref{luck}, so we will formulate it as a
proposition. The reader may however postpone the reading of the rest
of this section until then. Here we  mention related work of M.
Brion ~\cite{br} which exhibits degenerations of the diagonal in a
$G/P\times G/P$ into a union of products of Schubert varieties (as
below).

Let $K=\spane\{k_2,\dots,k_n\}$.   Define complete flags on $W$
$$F_{\bull}:0\subsetneq F_1\subsetneq F_2\subsetneq\dots\subsetneq F_n=W$$
$$F'_{\bull}:0\subsetneq F'_1\subsetneq F'_2\subsetneq\dots\subsetneq F'_n=W$$
 as
follows $F_1=L,\  F_{i}=L\ \oplus\ \spane\{k_2,\dots,k_{i}\}$ for
$i=2,\dots,n$, and $F'_i=\spane\{k_{n},\dots,k_{n-i+1}\}$ for
$i=1,\dots,n-1$ and $F'_n=W$. Define $\phi_{a}(t):W\to W$ for
$a=1,\dots,n$ as follows
\begin{itemize}
\item $\phi_1(t)$ is multiplication by $t$ on $L$ and identity on
$K$.
\item For $a=2,\dots,n$, $\phi_a(t)$ is multiplication by $t$ on $F_{a}$ and
$\phi_{a}(t)(k_i)=k_i$ for $i>a$.
\end{itemize}

\begin{proposition}\label{bord}
The following limits hold (the limit operations below are in the
Hilbert scheme of subschemes of $\Gr(r,W)\times \Gr(r,W)$). In the
formulas below $I$ varies over all subsets of $[n]$ of cardinality
$r$ (with additional conditions in Equations ~\eqref{twoo} and
~\eqref{three})
\begin{equation}\label{one}
\lim_{t_n\to 0} (1\times \phi_n(t_n)) \lim_{t_{n-1}\to 0} (1\times
\phi_{n-1}(t_{n-1}))\dots \lim_{t_1\to 0} (1\times
\phi_1(t_1))\Delta\ =\
\bigcup_{I}\Omega_I(F_{\bull})\times\Omega_{I'}(F'_{\bull})
\end{equation}
\begin{equation}\label{twoo}
\lim_{t_n\to 0} (1\times \phi_n(t_n)) \lim_{t_{n-1}\to 0} (1\times
\phi_{n-1}(t_{n-1}))\dots \lim_{t_2\to 0} (1\times
\phi_2(t_2))\deltaplus\ =\ \bigcup_{I:1\in
I}\Omega_I(F_{\bull})\times\Omega_{I'}(F'_{\bull})
\end{equation}
\begin{equation}\label{three}
\lim_{t_n\to 0} (1\times \phi_n(t_n)) \lim_{t_{n-1}\to 0} (1\times
\phi_{n-1}(t_{n-1}))\dots \lim_{t_2\to 0} (1\times
\phi_2(t_2))\deltaminus\ =\ \bigcup_{I:1\not\in
I}\Omega_I(F_{\bull})\times\Omega_{I'}(F'_{\bull})
\end{equation}
\end{proposition}
\begin{proof}
To prove this by induction, assume $W$ is an $n$ dimensional space,
$W=L\oplus K$ with $K$ $(n-1)$-dimensional. First degenerate
$\Delta\in \Gr(r,W)\times \Gr(r,W)$ to $\deltaplus\cup\deltaminus$.
 $\deltaplus$ (and $\deltaminus$) is a fiber bundle over the usual diagonal
 $\Delta(r-1,K)\subseteq
 \Gr(r-1,K)\times\Gr(r-1,K)$. We now degenerate $\Delta(r-1,K)$ (and
 similarly $\deltaminus$). The conclusions are now easy to prove.
\end{proof}
\subsection{Tangent spaces}

Let $$\deltaplus^o=\deltaplus-\deltaminus.$$ It is easy to see that
$\deltaplus^o$ is a smooth open dense subset of $\deltaplus$.

 Let $(A,B)\in \deltaplus $. Recall that $\beta:W\to W$ is the projection to $K$ corresponding to the decomposition
  $W=L\oplus K$. Then,
$\beta(A)\subset B$ and $\beta$ induces a morphism $W/A\to W/B$.
\begin{lemma}\label{maryam} For $(A,B)\in \deltaplus^o$, the induced map  $\beta:W/A\to W/B$ is an
isomorphism. The induced map $\beta:A\to B$ is singular with kernel
of rank $1$.
\end{lemma}
\begin{proof} Consider the natural map  $p:K\to W\to W/B$ (inclusion followed by projection). The rank of
the kernel is $\rk(K\cap B)=\rk(B)-1$ (since
$(A,B)\not\in\deltaminus$, $B\not\subseteq K$), and hence the rank
of the image is $\rk(K)+1-\rk(B)=\rk(W/B)$. Therefore the composite
$K\to W/B$ is surjective. But $\beta$ induces a surjective map $W\to
K$. The composite $W\leto{\beta} K\leto{p} W/B$ is surjective and A
is in its kernel. So $\beta$ induces a surjective (hence
isomorphism) map $W/A\to W/B$.

Since $\ker(\beta)=L$ and $L\subset A$ the second assertion is
clear.
\end{proof}
 \begin{lemma}\label{tangentlemma}
 The tangent space to the scheme $\Delta_{\beta}$ at $(A,B)$ is the
 vector subspace of the tangent space of $\Gr(r,W) \times \Gr(r,W)$ at $(A,B)$ given
by pairs of maps $A\leto{\Gamma_A} W/A,\ B\leto{\Gamma_B} W/B$ such
that the following diagram commutes:
$$\xymatrix{A\ar[r]^{\beta}\ar[d]^{\Gamma_A}  &                B\ar[d]^{\Gamma_B}\\
            W/A\ar[r]^{\beta}            &                    W/B }$$
\end{lemma}
\begin{proof}
 Choose splittings $W=A\oplus W/A$ and $W=B\oplus W/B$. For
$\psi\in W/A$, write $\beta(0\oplus\psi)=b(\psi)\oplus \beta(\psi)$
(the last term arises from the natural map $W/A\to W/B$).

 Suppose $\Gamma_A:A\to W/A$ and $\Gamma_B:B\to W/B$ are
deformations of $A$ and $B$ respectively. This means that the
deformed $A$ (similarly $B$) is the $\Bbb{C}[\epsilon]/(\epsilon^2)$
span of elements of the form $(a+\epsilon\Gamma_A(a))$. The
condition that this deformation stays inside $\Delta_{\beta}$ is
that $\beta(a+\epsilon\Gamma_A(a))$ should be in the
$\Bbb{C}[\epsilon]/(\epsilon^2)$ span of terms of  the form
$b+\epsilon\Gamma_B(b)$. A typical element in this span looks like
$b+\epsilon(b'+\Gamma_B(b))$. Write
$\beta(a+\epsilon\Gamma_A(a))=b+\epsilon(b'+\Gamma_B(b)).$
 This forces $b=\beta(a)$, and hence reading the $\epsilon$ term and its component in
 $W/B$,
we get the desired commutativity.
\end{proof}

\section{S-bundles and the diagonal}\label{beethoven}
Let $T$ be a ES-vector bundle of rank $r+k$ and degree $k(g-1)$ on
$\pone$. Fix isomorphisms $G_j:T_{p_j}\to W$ and  $H_j: T_{q_j}\to
W$.

We therefore have maps
$$\Gr(0,r,T)\to \prod_{j=1}^g\bigl(\Gr(r,T_{p_j})\times\Gr(r,T_{q_j})\bigr)\to
\prod_{j=1}^g\bigl(\Gr(r,W)\times\Gr(r,W)\bigr)$$

Call the composite $\Phi=\Phi(G,H)$. Fix $W=L\oplus K$ and the
projection $\beta:W\to W$ as in Section ~\ref{godfather}. Define
$\gamma_j:T_{p_j}\to T_{q_j}$ using the diagram
\begin{equation}\label{wed}
\xymatrix{T_{p_j}\ar[r]^{G_j}\ar[d]^{\gamma_j}      &    W\ar[d]^{\beta} \\
               T_{q_j}\ar[r]^{H_j}              &                                  W\\
            }
            \end{equation}

\begin{proposition}\label{shemoved}
\begin{enumerate}
\item[(a)] For generic $(G,H)$, $\Phi(G,H)^{-1}(\prod_{j=1}^g\deltaplus)=\Phi(G,H)^{-1}(\prod_{j=1}^g\deltaplus^o).$
\item[(b)] $\Phi(G,H)^{-1}(\prod_{j=1}^g\deltaplus)$ is in bijection
with the set of singular subbundles of degree $0$ and rank $r$ of
$(T,\gamma)$ considered in Section ~\ref{before}.
\item[(c)] For generic $(G,H)$, $\Phi(G,H)^{-1}(\prod_{j=1}^g\deltaplus^o)$ is
reduced and of the expected dimension $0$.
\end{enumerate}
\end{proposition}
\begin{proof}
Item (b) follows easily from the definitions.

The expected dimension is $$\dim\Gr(0,r,T) - g\dim\Gr(r,W)\ =\
rk+k(g-1)r-grk=0.$$

Note that $\GL(W)\times\GL(W)$ acts transitively on
$\Gr(r,W)\times\Gr(W)$. Therefore the group
$\prod_{j=1}^g\bigl(\GL(W)\times\GL(W)\bigr)$ acts transitively on
$\prod_{j=1}^g\bigl(\Gr(r,W)\times\Gr(r,W)\bigr)$. We now recall
that $\deltaplus^o$ is smooth and connected. Hence, by Kleiman's
transversality theorem (cf. ~\cite{kl}, ~\cite{int} \S B.9.2), for
generic $\tau\in\prod_{j=1}^g\bigl(\GL(W)\times\GL(W)\bigr)$,
$\Phi^{-1}\tau^{-1}(\prod_{j=1}^g\deltaplus^o)$ is reduced and of
the expected dimension $=0$. Notice that the action of $\tau$ just
modifies the maps $G$ and $H$. Therefore the assertion (c) follows.
The equality in (a) holds because if $S$ is a nonempty subset of
$\{1,\dots,g\}$,  the expected dimension of (any component of)
$\Phi(G,H)^{-1}(\prod_{j=1}^g A_j)$ is negative where
$A_j=\deltaplus$ if $j\not \in S$ and $A_j=\deltaplus
\cap\deltaminus$ if $j\in S$.
\end{proof}
 Now let $G$ and $H$ be generic and
 $$m=|\Phi^{-1}(\prod_{j=1}^g \deltaplus^o)|,\ \{V^{(1)},\dots,V^{(m)}\}=\Phi^{-1}(\prod_{j=1}^g \deltaplus^o).$$
Let $\mqq^{(a)}=T/V^{(a)}$ for $a=1,\dots,m$. We will now endow
$V^{(a)}$ and $\mqq^{(a)}$ with S-bundle structure using the
diagrams, obtaining $S$-bundle $(V^{(a)},\tau^{(a)})$ and
$(Q^{(a)},\theta^{(a)})$ respectively:
\begin{equation}\label{one}
\xymatrix{V^{(a)}_{p_j}\ar[r]\ar[d]^{\tau^a_{j}}  &     T_{p_j}\ar[r]^{G_j}\ar[d]^{\gamma_{j}}      &    W\ar[d]^{\beta} \\
            V^{(a)}_{q_j}\ar[r]                       &   T_{q_j}\ar[r]^{H_j}              &
            W\\
            }
            \end{equation}

            and
\begin{equation}\label{two}
\xymatrix{T_{p_j}\ar[r]\ar[d]^{\gamma_{j}}  &     \mqq^{(a)}_{p_j}\ar[d]^{\theta^{(a)}_{j}} \\
            T_{q_j}\ar[r]                     &     \mqq^{(a)}_{q_j}\\
            }
\end{equation}

\begin{proposition}\label{hotri}
The inequality $m\geq M(r,k,g)$ holds.
\end{proposition}
The proof of Proposition ~\ref{hotri} will appear in Section
~\ref{agni}. It will eventually be shown $m=M(r,k,g)$.

\begin{lemma}\label{tangentcalcul}
The tangent space at $V^{(a)}$ to
$\Phi^{-1}(\prod_{j=1}^g\deltaplus^o )$ is
$\uhom((V^{(a)},\tau^{(a)}),(\mqq^{(a)},\theta^{(a)}))$.
\end{lemma}
\begin{proof}
The tangent space to $\Gr(0,r,T)$ at $V^{(a)}\subseteq T$ is
$\home(V^{(a)},T/V^{(a)})$ by Grothendieck's theory of quot schemes.

For $j=1,\dots,g$, let $\Phi_j$ be the composition (where the last
map  is the projection to the $j$th factor):
$$\Gr(0,r,T)\to
\prod_{j=1}^g\Gr(r,T_{p_j})\times\Gr(r,T_{q_j})\leto{\prod(G_j,H_j)}
\prod_{j=1}^g\Gr(r,W)\times\Gr(r,W)\to \Gr(r,W)\times\Gr(r,W).$$
 An element $\Gamma \in\home(V^{(a)},T/V^{(a)})$ is in the tangent
space of $\Phi^{-1}(\prod_{j=1}^g\deltaplus )$ if and only for each
$j=1,\dots,g$, $(\Phi_j)_*(\Gamma)$ is in the tangent space of
$\deltaplus^o$ at $(G_j(V^{(a)}_{p_j}),H_j(V^{(a)}_{q_j})).$ By
Lemma ~\ref{tangentlemma}, this condition implies that, as desired,
$\Gamma\in\uhom((V^{(a)},\tau^{(a)}),(\mqq^{(a)},\theta^{(a)}))$.
\end{proof}
\begin{lemma}\label{versace}
\begin{enumerate}
\item[(i)] $\uhom((V^{(a)},\tau^{(a)}),(\mqq^{(b)},\theta^{(b)}))\neq 0$
if and only if $a\neq b$.
\item[(ii)] The coherent sheaves underlying the  $S$-bundles $(\mqq^{(a)},\theta^{(a)}),\ a=1,\dots,m$ are vector
bundles on  $C$.
\item[(iii)] For $a=1,\dots, m$ and $j=1,\dots,g$, the linear map of vector spaces $\tau^{(a)}_{j}: V^{(a)}_{p_j}\to
V^{(a)}_{q_j}$ is singular with a rank $1$ kernel.
\item[(iv)] The vector bundles $V^{(a)}$ (resp. $Q^{(a)}$) for $a=1,\dots,m$ on
$\pone$ are ES and are hence isomorphic. In particular,
$V^{(a)}\leto{\sim}\mathcal{O}_{\pone}^{\oplus r}$.
\end{enumerate}
\end{lemma}
\begin{proof} (The proof models a similar proof in
~\cite{invariant}) By Lemma ~\ref{tangentcalcul},  the tangent space
at $V^{(a)}$ to $\Phi^{-1}(\deltaplus )$  is
$\uhom((V^{(a)},\tau^{(a)}),(\mqq^{(a)},\theta^{(a)}))$
 which is consequently $0$, because of Proposition ~\ref{shemoved}.

If $a\neq b$ consider the nonzero (!) composite $V^{(a)}\to
T\to\mqq^{(b)}$ (inclusion followed by projection). This composite
belongs to $\uhom((V^{(a)},\tau^{(a)}),(\mqq^{(b)},\theta^{(b)}))$
(the property of being a morphism of S-bundle follows from Diagrams
~\eqref{one} and ~\eqref{two}). This proves (i).

Now $\mvv^{(a)}\in \Phi^{-1}(\deltaplus^o)$, and therefore
 $\theta^{(a)}_{j}:\mqq^{(a)}_{p_j}\to Q^{(a)}_{q_j}$ is an isomorphism
(see Lemma ~\ref{maryam}). So (ii) follows immediately.

Assertion  (iii), follows because $L\subseteq\ker(\beta)$. Assertion
(iv) follows from Proposition ~\ref{bigone} and Kleiman's
transversality theorem.
\end{proof}

\section{S-bundles and determinant of cohomology}\label{hin}
Let $V=\mathcal{O}_{\pone}^r$ and $Q$ an ES bundle of degree
$k(g-1)$ and rank $k$.  Let $(V,\tau)$ and $(Q,\theta)$ be S-bundles
on $C$ so that $\theta_j$ is an  isomorphism for
 $j=1,\dots,g$. It follows that the coherent sheaf underlying $(Q,\theta)$ is a honest vector bundle
$\tilde{Q}$ on $C$.

 Define $\ma=\mathcal{A}(\tau,\theta)$ a coherent sheaf on $C$
 by the following exact sequence of sheaves:
 $$0\to \mathcal A\to f_{*}\shom(V,Q)\leto{\Sigma} \oplus_{j=1}^g \home(V_{p_j},Q_{q_j})|r_j\to 0$$
where the last sheaf is a skyscraper sheaf at $r_j$ and
$$\Sigma
(\Gamma)=\theta_j\circ\Gamma_{p_j}- \Gamma_{q_j}\circ\tau_{j}\in
\home(V_{p_j},Q_{p_j}).$$

Since $\theta_j$ are assumed to be isomorphisms, the map $\Sigma$ is
 surjective. If $\tau$ and $\theta$ vary in families so that
$\theta_j$ are isomorphisms (on each fiber over the parameter space)
then $\ma(\tau,\theta)$ is flat over the parameter space. This is
because the kernel of a surjective map of flat modules is flat (long
exact sequence in Tor!).

Now assume that $\deg(V)=0$ and $\deg(Q)=k(g-1)$ then since $f$ is a
finite morphism,
$$\chi(C,\ma)=\chi(C,f_* \shom(V,Q)) -grk
=\chi(\pone,\shom(V,Q))-rk= kr(g-1) + kr-grk =0.$$ So there is then
a canonical section
\begin{equation}\label{cano}\sigma(\tau,\theta)\in
\D(\ma(\tau,\theta))\leto{\sim}\D(f_*(\shom(V,Q)))\tensor
\tensor_{j=1}^g\det(\home(V_{p_j},Q_{q_j}))^*.
\end{equation}
We note the following properties of this construction.
 \begin{lemma}\label{jaan}
\begin{enumerate}
\item $H^0(C,\ma)=\home((V,\tau),(Q,\theta))$.
\item $\sigma(\tau,\theta)\neq 0$ if and only if
$\home((V,\tau),(Q,\theta))= 0$.
\item If the $\tau_j$ are
isomorphisms for $j=1,\dots,g$, and $\tilde{V}$ the locally free
coherent sheaf underlying $(V,\tau)$ then we have a natural
isomorphism of sheaves on $C$: $\ma\to \shom(\tilde{V},\tilde{Q})$.
\end{enumerate}
\end{lemma}

\subsection{Geometric S-bundles}\label{agree}
Now consider the S-bundles obtained from geometry in Section
~\ref{before} (Lemma ~\ref{versace}). By Lemma ~\ref{versace}, (iv),
the vector bundles $V^{(a)}$ on $\pone$ are all isomorphic to the
vector bundle $V=\mathcal{O}^{\oplus r}$, and the bundles
$\mqq^{(a)}$ are all isomorphic to $Q$ where $Q$ is the unique ES
bundle of the degree and rank of $\mqq^{(a)}$. We choose
isomorphisms
 $(V^{(a)},\tau^{(a)})\leto{\sim} (V,\tilde{\tau}^{(a)})$ and
$(Q^{(a)},\theta^{(a)})\leto{\sim}(Q,\tilde{\theta}^{(a)})$ for
$a=1,\dots,m$; for suitable $\tilde{\tau}$ and $\tilde{\theta}$.

From Lemma ~\ref{versace} we know that the  $\ttheta^{(b)}_j$ are
isomorphisms for $j=1,\dots,g$. Consider,
$$\sigma(\ttau^{(a)},\ttheta^{(b)})\in \D(\ma(\ttau^{(a)},\ttheta^{(b)}))\leto{\sim}\D(f_*(\shom(V,Q)))\tensor
\tensor_{j=1}^g\det(\home(V_{p_j},Q_{q_j}))^*.$$

We may think of this as a $m\times m$ matrix with entries in the
 one dimensional vector space $\D(f_*(\shom(V,Q)))\tensor
\tensor_{j=1}^g\det(\home(V_{p_j},Q_{q_j}))^*$. By Lemma
~\ref{versace} and Lemma ~\ref{jaan} (2), the $(a,b)$ entry is
non-zero exactly when $a=b$. Hence the matrix is non-singular.

Find maps ${\lambda}^{(a)}_{j}(t)\in \home(V_{p_j},
V_{q_j})=\home(\Bbb{C}^r,\Bbb{C}^r)$ (using the identification $V\to
\mathcal{O}^{\oplus r}$) so that for $j=1,\dots,g$,
\begin{enumerate}
\item $\lambda^{(a)}_{j}(0)=\ttau^{(a)}_{j}$
\item $\det(\lambda^{(a)}_{j}(t))= t$.
\end{enumerate}
 (we may find such matrices using Jordan
canonical forms because the kernel of $\ttau^{(a)}_j$ is one
dimensional for $a=1,\dots,m$ and $j=1,\dots,g$). Therefore, for
$a,a'\in\{1,\dots,m\}$,
$\det(\lambda^{(a)}_{j}(t))=\det(\lambda^{(a')}_{j}(t))$ and these
numbers are non-zero if $t\neq 0$. Consider the $m\times m$ matrix
formed by
$$\sigma(\lambda^{(a)}(t),\ttheta^{(b)})\in \D(\ma(\lambda^{(a)}(t),\ttheta^{(b)}))\leto{\sim}\D(f_*(\shom(V,Q)))\tensor
\tensor_{j=1}^g\det(\home(V_{p_j},Q_{q_j}))^*.$$

This matrix is going to be non-singular for values of $t$ in a
sufficiently small Zariski neighborhood of  $t=0$. Let $\epsilon\neq
0$ be one such value for $t$. Let $\tv^{(a)}$ be the locally free
coherent sheaf underlying $(V,\lambda^{(a)}(\epsilon))$. Let
$\tq^{(b)}$ be coherent sheaf on $C$ underlying $(Q,\ttheta^{(b)})$.
From Lemma ~\ref{versace}, we see that $\tq^{(b)}$ are locally free.
According to Lemma ~\ref{looks}, $\tv^{(a)}$ have isomorphic
determinant line bundles of degree $0$ for $a=1,\dots,m$.
\begin{lemma}
The tuple
$((\tv^{(1)})^*,\dots,(\tv^{(m)})^*;\tq^{(1)},\dots,\tq^{(m)})$ is
an $(m,r,k)$-frame on the  rational nodal curve $C$.
\end{lemma}
Together with Lemma ~\ref{moment} and Proposition ~\ref{map}, this
would conclude the proof of the strange duality conjecture for
generic curves if we can show that the identifications that we made
here are compatible with the determinant operation that went with
the definition of the frame.
\subsection{The basic compatibility verification}
Let $(V,\tau)$ and $(V,\eta)$ be $S$-bundles such that $V$ is of
degree $0$ and rank $r$. Let $\tilde{V}^1$ and $\tilde{V}_2$ be the
underlying coherent sheaves on $C$. Assume that $\tau_j$ and
$\eta_j$ are isomorphisms for $j=1,\dots,g$ (so $\tilde{V}^1$ and
$\tilde{V}^2$ are vector bundles on $C$). Assume further that
$\det(\tilde{V}^1)\leto{\sim}\det(\tilde{V}^2)$. Let $(Q,\theta)$
and $(Q,\delta)$ be $S$-bundles such that $\theta_j$ and $\delta_j$
are isomorphisms $j=1,\dots,g$ and $\tq_1$ and $\tq_2$ the
corresponding vector bundles on $C$. Assume that $Q$ is of rank $k$
and degree $k(g-1)$. We need to verify that the following diagram
commutes:
\begin{equation}\label{correcto}
\xymatrix@R=0.05cm{\D(\tv^*_1\tensor\tq_1)\tensor\D(\tv^*_2\tensor\tq_2) \ar[ddd]\ar[r]^{\Phi_{\tv^*_1,\tv^*_2}} & \D(\tv^*_2\tensor\tq_1)\tensor\D(\tv^*_1\tensor\tq_2)\ar[ddd]\\
\\
\\
\D(f_*(\shom(V, Q)))\tensor\prod_{j=1}^g\det(\home(V_{p_j},
Q_{q_j}))^*\ar[r]^{\operatorname{=}} &
\D(f_*(\shom(V,Q)))\tensor\prod_{j=1}^g\det(\home(V_{p_j},
Q_{q_j}))^*\\
\text{   }\tensor   & \text{   }\tensor\\
\D(f_*(\shom(V, Q)))\tensor\prod_{j=1}^g\det(\home(V_{p_j},
Q_{q_j}))^* &
\D(f_*(\shom(V,Q)))\tensor\prod_{j=1}^g\det(\home(V_{p_j},
Q_{q_j}))^*}
\end{equation}

For this we first consider two natural exact sequences associated to
$\tq_1$ and $\tq_2$. The one associated to $\tq_1$ is
\begin{equation}
0\to \tq_1\to f_*Q \leto{\alpha}\oplus_{j=1}^g Q_{q_j}\mid r_j\to 0
\end{equation}
where $\alpha_j:Q_{p_j}\oplus Q_{q_j}\to Q_{q_j}$ takes
$(q_1,q_2)\to q_2-\theta_{j}(q_1)$. Notice that for $a=1,2$ there
are isomorphisms $f_*Q\tensor \tv^*_a\to f_*(V^*\tensor
Q)=f_*\shom(V,Q)$, and $(Q_{q_j}\mid r_j)\tensor\tv^*_a=
Q_{q_j}\tensor V^*_{p_j}\mid r_j=\home(V_{p_j},Q_{q_j})\mid r_j$.
From Section ~\ref{falt}, recall the second observation, and the
compatibility under exact sequences in $(I,J)$. We now see that
diagram ~\eqref{correcto} commutes.
\section{Proof of
Proposition ~\ref{hotri}}\label{agni}
\subsection{Conformal blocks}
Irreducible polynomial representations of  $U(r)$ are parameterized
by weakly decreasing sequences of non-negative integers
$\lambda=(\lambda_1\geq\lambda_2\geq\dots\geq\lambda_{r})\in\Bbb{Z}^r$.
These restrict to irreducible representations $\bar{\lambda}$ of
$\SU(r)$. Also note that$(\lambda_1,\lambda_2,\dots,\lambda_r)$ and
$(\mu_1,\mu_2,\dots,\mu_r)$ restrict to give the same irreducible
representation of $\SU(r)$ if $\lambda_a-\mu_a$ is a constant for
$a=1,\dots,r$.  The congruence class of $|\bar{\mu}|=\sum_a\mu_a$
$\pmod{r}$ is therefore a well defined
 $\Bbb{Z}/r\Bbb{Z}$ ``invariant'' of the representation $\bar{\mu}$ of $\SU(r)$.

The dual of a representation $\bar{\mu}$ of $\SU(r)$ is denoted by
$\bar{\mu}^*$ and equals $(\mu_1-\mu_r, \mu_1-\mu_{r-1}, \dots,
\mu_1-\mu_2, 0)$.

A representation $\bar{\lambda}$ of $\SU(r)$ with
$\lambda=(\lambda_1,\lambda_2,\dots,\lambda_{r})$ is said to have
level $\leq k$ if $\lambda_1-\lambda_r\leq k$.

Given irreducible representations $\bar{\mu}^{1},\dots,
\bar{\mu}^{s}$ of level $\leq k$ of  $\SU(r)$,  we obtain the
dimensions of spaces of conformal blocks
$N^{(k)}_g(\bar{\mu}^{1},\dots,\bar{\mu}^{s})$ (as in
~\cite{beauu}).

It is known that $N^{(k)}_g(\emptyset)=M(r,k,g)$(cf.
~\cite{v2},~\cite{v1} and ~\cite{v3}).

The factorization theorem of Tsuchiya-Ueno-Yamada (see ~\cite{TUY}
and ~\cite{beauu})  gives
$$N^{(k)}_g(\bar{\mu}^{1},\dots,\bar{\mu}^{s})=\sum_{\bar{\nu}} N^{(k)}_{g-1}(\bar{\mu}^{1},\dots,\bar{\mu}^{s},\bar{\nu},\bar{\nu}^*)$$
(where $\bar{\nu}$ runs through all irreducible representations of
$\SU(r)$ of level $\leq k$.) Repeated application of this formula
gives
\begin{equation}\label{tuytuy}
N^{(k)}_g(\emptyset)=\sum_{\bar{\mu}^{1},\dots,\bar{\mu}^{g}}
N^{(k)}_0(\bar{\mu}^{1},(\bar{\mu}^{1})^*,\bar{\mu}^{2},(\bar{\mu}^{2})^*,\dots,\bar{\mu}^{g},(\bar{\mu}^{g})^*)
\end{equation}
where $\bar{\mu}^{1},\dots,\bar{\mu}^{g}$ vary over all irreducible
representations of $\SU(r)$ of level $\leq k$.

The Verlinde algebra for $\SU(r)$ at level $\leq k$ is given by the
rule
$$\bar{\mu}^{1}\cdot\bar{\mu}^{2}\cdots \bar{\mu}^{s}=\sum _{\bar{\nu}}
N^{(k)}_0(\bar{\mu}^{1},\dots,\bar{\mu}^{s},\bar{\nu})\bar{\nu}^*$$
where $\bar{\nu}$ runs through all irreducible representations of
$\SU(r)$ of level $\leq k$.
\subsection{Generalized Gromov-Witten numbers}\label{temina}
Fix positive integers $r$ and $k$  and set $n=r+k$. Suppose
$I^{1},\dots,I^{s}$ are  subsets of $[n]=\{1,\dots,n\}$ each of
cardinality $r$.

Let ${\pto}_1,\dots,{\pto}_s$ be a set of distinct points on
$\pone$. Let $T$ be an evenly split (ES) bundle of degree $-D$ and
rank $n$ (recall the definition of ES bundles given in Section
~\ref{before}). Choose generic complete flags $E^{{\pto}_j}_{\bull}$
on the fibers $T_{{\pto}_j}$ for $j=1,\dots,s$.

Define
$\langle\omega_{I^1},\omega_{I^2},\dots,\omega_{I^s}\rangle_{d,D}$
to be the number of subbundles ($0$ if this number is infinite) $V$
of $T$ of degree $-d$ and rank $r$, such that for $j=1,\dots,s$, the
point $V_{{\pto}_j}\in \Gr(r,T_{{\pto}_j})$ lies in the Schubert
cell $\Omega^o_{I^{j}}(E^{{\pto}_j}_{\bull})$. By Kleiman's
transversality theorem the space of such $V$ is equidimensional of
dimension (using Proposition  ~\ref{bigone})
$$\dim(\Gr(d,r,T))-\sum_{j=1}^s\codim(\omega_{I^j})=r(n-r)+dn-Dr-\sum_{j=1}^s\codim(\omega_{I^j}).$$

 If $D=0$ it is easy to see that the above definition gives the structure coefficients in the small
 quantum cohomology of $\Gr(r,n)$ (for example see ~\cite{FuP}). Here we use the standard bijection between  subbundles of $\mathcal{O}_{\pone}^{\oplus n}$
 of rank $r$ and degree $-d$ and  maps $\pone\to
 \Gr(r,n)$ of degree $d$. The numbers when $D\neq 0$ can be
 recovered from the small quantum cohomology structure constants by
 using shift operations (see Proposition ~\ref{shiftt1}).

 One also notes that if
 $L=\{n-r+1,n-r+2,\dots,n\}$ then
 \begin{equation}\label{vacuum}
\langle\omega_{I^1},\omega_{I^2},\dots,\omega_{I^s}\rangle_{d,D}\ =\
\langle\omega_{I^1},\omega_{I^2},\dots,\omega_{I^s},\omega_L\rangle_{d,D}
\end{equation}
 This is because we  are  imposing
an open dense condition at the $s+1^{\operatorname{th}}$ point.

\subsection{Proof of Proposition ~\ref{hotri}}\label{luck}
We return to the notation and setting of Proposition ~\ref{hotri}.
We introduce a new piece of notation relating Schubert cells in
$\Gr(r,n)$ and representations of $\GL(r)$. To a $r$-element subset
$I=\{i_1<\dots<i_r\}$ of $[n]=\{1,\dots,n\}$ we associate a weakly
decreasing sequence of non-negative integers
\begin{equation}\label{convent}
I\mapsto
\lambda(I)=(\lambda_1\geq\lambda_2\geq\dots\geq\lambda_r)\in\Bbb{Z}_{\geq
0}^r,
\end{equation}
where $\lambda_a=n-r+a-i_a,\ a=1,\dots,r.$ We note the obvious:
\begin{lemma}\label{mallika}
\begin{enumerate}
\item  As $I$ runs through all $r$-element subsets of $[n]$ with $1\in I$,
$\bar{\lambda}(I)$ runs (exactly once) through all irreducible
representations of $\SU(r)$ of level $\leq k=n-r$.
\item For $I$ as in (1),
$\bar{\lambda}(I')= \bar{\lambda}(I)^*.$
\end{enumerate}
\end{lemma}

We return to the Proof of Proposition ~\ref{hotri}. Recall that by
Proposition ~\ref{bord}, $\deltaplus $ degenerates into
$$\bigcup_{I:i_1=1}\Omega_I(F_{\bull})\times\Omega_{I'}(F'_{\bull})\subseteq
\Gr(r,W)\times\Gr(r,W).$$ We may assume by modifying the maps $G_j$,
$H_j$ in diagram ~\eqref{wed} that any
$\Phi^{-1}(\prod_{j=1}^g\Omega_{I^j}(F_{\bull})\times\Omega_{(I^j)'}(F'_{\bull}))$
is smooth of the expected dimension zero (using Kleiman's
transversality theorem). We may assume that any point of
intersection is in
 $\Phi^{-1}(\prod_{j=1}^g\Omega^o_{I^j}(F_{\bull})\times\Omega^o_{(I^j)'}(F'_{\bull})).$

 But this number is just
$$\langle\omega_{I^1},\dots,\omega_{I^s},\omega_{(I^1)'},\dots,\omega_{(I^g)'}\rangle_{0,-k(g-1)}$$

This implies that for some $(t_1,\dots,t_n)$ with $t_i\neq 0$,
$$(1\times \phi_n(t_n))(1\times
\phi_{n-1}(t_{n-1}))\dots (1\times \phi_1(t_1))\deltaplus$$ has at
least
$$m_0=\sum\langle\omega_{I^1},\dots,\omega_{I^g},\omega_{(I^1)'},\dots,\omega_{(I^g)'}\rangle_{0,-k(g-1)}$$
transverse points of intersections, where the sum is over all
sequences  of subsets $(I^1,\dots,I^g)$ of $[n]$ with $r$ elements,
each of which contain $1$ (see equation ~\eqref{twoo}).

Hence for generic $G$ and $H$ (by absorbing $1$ and
$\phi_n(t_n)\circ\dots\phi(t_1)$ in the maps $G$ and $H$), we see
that $\Phi(G,H)^{-1}(\prod_{j=1}^g\deltaplus^o)$ has at least  $m_0$
transverse points (the passage from $\deltaplus$ to $\deltaplus^o$
is clear because the intersection is over smooth points of
$\deltaplus$). So $m\geq m_0$.

 The following Proposition follows from a  theorem of
Witten  and will be proved in  Section ~\ref{witten}.
\begin{proposition}\label{hoopdreams}
Suppose  $I^1,\dots,I^g$ are $r$-element subsets of $[n]$, and
$\mu^j=\lambda(I^j)$ for $j=1,\dots,g$. Then,
$$\langle\omega_{I^1},\dots,\omega_{I^s},\omega_{(I^1)'},\dots,\omega_{(I^g)'}\rangle_{0,-k(g-1)}=
N^{(k)}_0(\bar{\mu}^{1},(\bar{\mu}^{1})^*,\bar{\mu}^{2},(\bar{\mu}^{2})^*,\dots,\bar{\mu}^{g},(\bar{\mu}^{g})^*)$$
\end{proposition}
From Proposition ~\ref{hoopdreams} and Lemma ~\ref{mallika} (1), it
follows that
$$m_0=\sum_{\bar{\mu}^{1},\dots,\bar{\mu}^{g}}
N^{(k)}_0(\bar{\mu}^{1},(\bar{\mu}^{1})^*,\bar{\mu}^{2},(\bar{\mu}^{2})^*,\dots,\bar{\mu}^{g},(\bar{\mu}^{g})^*)$$
where $\bar{\mu}^{1},\dots,\bar{\mu}^{g}$ vary over all irreducible
representations of $\SU(r)$ of level $\leq k$. By Equation
~\eqref{tuytuy}, one sees now that $m_0= M(r,k,g)$ and this
concludes the proof of Proposition ~\ref{hotri}.
\subsection{Examples}\label{examples}
We first consider the case $k=1$. In this case in the formula
$(\ddagger)$ $n=r+1$ and $I^1,\dots, I^g$ are subsets of
$[r+1]=\{1,\dots,r+1\}$ each with $r$ elements and containing $1$.
There are $r$ choices for each $I^j$. It is an easy calculation
using the shift operations from Proposition ~\ref{shiftt1} that each
of the summands of $(\ddagger)$ is $1$ and hence we obtain sum $r^g$
as expected.

We move on to a more non trivial example. We compute the sum in
$(\ddagger)$ to see that it agrees with $2^{g-1}(2^g+1)$ (see
~\cite{beauu}, Section 9). Here we notice that there are three
allowable $I$'s -$\{1,2\}$, $\{1,3\}$ and $\{1,4\}$. By a small
calculation we obtain: If a summand in $(\ddagger)$ has $\ell$
$\{1,3\}$'s , then that summand is $1$ if $\ell=1$ or $2$ and
$2^{\ell-1}$ otherwise. So the sum  is
$$\left( \begin{array}{cc} g \\ 0
\end{array}\right) 2^g +\left( \begin{array}{cr} g \\ 1
\end{array}\right)2^{g-1}+\sum_{\ell=2}^g \left( \begin{array}{cc} g \\ \ell
\end{array}\right)2^{g-1}=2^{g-1}(2^g+1)$$

The number ${h^0(SU_X(r),\ml^k)}/{r^g}$ is known to be symmetric in
$k$ and $r$. It is possible to see this as a consequence of
Grassmann duality $\Gr(r,r+k)\leto{\sim}\Gr(k,r+k)$. The factor of
$r^g$ in the denominator comes from the requirement that
$I^1,\dots,I^g$ in $(\ddagger)$ each contain $1$ (This condition is
not symmetric under Grassmann duality - but symmetric in a cyclic
sense).
\subsection{Witten's theorem and consequences}
Recall that the small quantum cohomology of $\Gr(r,n)$ is an
associative ring  $\QH(\Gr(r,n))$ whose underlying abelian group is
$H^*(\Gr(r,n),\Bbb{Z})\tensor \Bbb{Z}[q]$ and the
 product structure is given by
  $$\omega_I\star\omega_J=\sum_{K,d}\langle
\omega_I,\omega_J,\omega_K\rangle_{d,0} \ q^d \omega_{K'}$$ (see
Section ~\ref{fultone} for the notation) where $K$ runs through all
subsets of $[n]$ of cardinality $r$, and $d$ runs though all non
negative integers and $\omega_{K'}$ is the dual of $\omega_K$, see
Section ~\ref{fultone}. We note the following generalization of the
 (small) quantum product of Schubert classes
$$\omega_{I^1}\star\omega_{I^2}\star\dots\star\omega_{I^s}=\sum_{K,d}\langle
\omega_{I^1},\omega_{I^2},\dots,\omega_{I^s},\omega_K\rangle_{d,0} \
q^d \omega_{K'}$$

Recall that $n=r+k$.  Witten's theorem ~\cite{witten} gives an
isomorphism
\begin{equation}\label{spice}
W:\operatorname{QH}^*(Gr(r,n))/(q-1)\to \Ri
\end{equation}
from the quantum cohomology of $\Gr(r,n)$ at $q=1$  to the Verlinde
algebra of  $U(r)$ at $SU(r)$ level $\leq k$ and $U(1)$ level $\leq
n$. $\Ri$ is additively generated by sequences
$\lambda=(\lambda_1\geq\lambda_2\geq\dots\geq\lambda_{r})\in\Bbb{Z}^r$
such that $\lambda_r\geq 0$ and $\lambda_1\leq k$. where
$\lambda_a=n-r+a-i_a,\ a=1,\dots,r.$ $W$ takes $\omega_I$ to the
partition $\lambda(I)$.

$R(U(r))_{k,n}$ is related to the product
 of the Verlinde algebra  $R(SU(r))_{k}$ for $\SU(r)$ at level $\leq k$
and the Verlinde algebra $R(U(1))_{rn}$ of $U(1)$ at level $\leq
nr$. This relation stems from the covering $\pi:\SU(r)\times U(1)\to
U(r)$. We will make this relation precise. $R(\SU(r))_k$ is a ring
with additive basis given by irreducible $\SU(r)$ representations of
level $\leq k$. $R(U(1))_{rn}$ is generated by $x$ with relation
$x^{nr}=1$.

Inside $R(SU(r))_{k}\times R(U(1))_{rn}$ consider the subspace
$\tilde{R}$ spanned by ${\bar{\lambda}}\times x^a$ so that  $a\equiv
|\bar{\lambda}|\pmod{r}$. This corresponds to representations of
$\SU(r)\times U(1)$ that descend to representations of $U(r)$ under
the covering $\pi$.  $\tilde{R}$ is a subring of $R(SU(r))_{k}\times
R(U(1))_{rn}$.

On $\tilde{R}$ (with the product structure) consider the  operator
\begin{equation}\label{word}
 T({\bar{\lambda}}\times x^a)= \bar{\eta}\times x^{n+a}
\end{equation}
 where
$\bar{\eta}$ is related to $\bar{\lambda}$  by the cyclic shift
$\eta=(k+\lambda_r\geq\lambda_1\geq\lambda_2\geq\dots\geq\lambda_{r-1})$
and $a\equiv |\bar{\lambda}|\pmod{r}$. It is easy to check that
$T(u)=u\cdot((k,0,\dots,0)\times x^n).$ The $\Bbb{Z}$-submodule of
$\tilde{R}$ generated by elements $T(u)-u$ for $u\in \tilde{R}$ is
an ideal $I$ of $\tilde{R}$ . From an easy calculation one can see
that each orbit $\{T^b(\bar{\mu}\times x^a)_{b=1}^r\}$ where
$a\equiv |\bar{\mu}|\pmod {r}$, contains a unique element of the
form $\bar{\lambda}\times x^{|\lambda|}\in \tilde{R}$ where the
sequence $\lambda$ satisfies $\lambda_1\leq k$ and $\lambda_r\geq 0$
( and $|\lambda|=\sum_{a=1}^r \lambda_a\in \Bbb{Z}$). The linear
$\tilde{R}\to \Ri$ sending $\bar{\lambda}\times x^{|\lambda|}\in
\tilde{R}$ to $\lambda\in \Ri$ is a ring homomorphism with kernel
$I$ and induces an isomorphism of rings $${\tilde{R}}/I \leto{\sim}
\Ri.$$

Let $K=\{1,k+2,k+3,\dots,k+r=n\}\subset[n]$. Then using Equation
~\eqref{word} on $\lambda=(0, 0,\dots, 0)\in \Bbb{Z}^r$,
\begin{equation}\label{domingo}
W(\omega_K)x^r=1\times 1
\end{equation}
\begin{remark} The reader  can verify the form of
Witten's theorem needed in this paper by just taking $\tilde{R}/I$
to the be the definition of $\Ri$, and deducing that the map $W$ is
a ring homomorphism from knowing that Pieri's rule holds on both
sides of Equation ~\eqref{spice} (using ~\cite{bert} on the quantum
cohomology side, and ~\cite{gepner} for the  Verlinde algebra for
$\SU(r)$ at level $\leq k$.)
\end{remark}
\subsubsection{Consequence of  Witten's theorem}\label{witten}
Let  $I^{1},\dots,I^{s}$ be subsets of $[n]$ each of cardinality
$r$. Assume that there is an integer $D$  so that
$\sum_{j=1}^s\codim(\omega_{I^s})=rk-Dr$. Then we have the following
Proposition which implies Proposition ~\ref{hoopdreams} immediately:
\begin{proposition}\label{sounds}
\begin{equation}\label{towit}
\langle\omega_{I^1},\omega_{I^2},\dots,\omega_{I^s}\rangle_{0,D} =
N^{(k)}_0(\bar{\mu}^{1},\dots,\bar{\mu}^{s})
\end{equation}
where $\bar{\mu}^{j}=\bar{\lambda}(I^j)$ are the associated
representations of $\SU(r)$.
\end{proposition}
\begin{proof}For simplicity assume $D\leq 0$ (this is the case when we need this
proposition). Let $K=\{1,k+2,k+3,\dots,n\}$. We claim
\begin{equation}\label{factor}
\langle\omega_{I^1},\omega_{I^2},\dots,\omega_{I^s}\rangle_{0,D}\ =
\
\langle\omega_{I^1},\omega_{I^2},\dots,\omega_{I^s},\omega_K,\omega_K,\dots,\omega_K\rangle_{-D,0}
\end{equation}
with $\omega_K$ repeated $-D$ times. To see this apply the shift
operation from Proposition ~\ref{shiftt1} to each of the $-D$
 $\omega_K$'s appearing on the right hand side of Equation
 ~\eqref{factor}.

 According to Witten's theorem the two sides of Equation
~\eqref{factor} equal the coefficient of $W(\omega_{\{1,\dots,r\}})$
in the product
\begin{equation}\label{realise}
W(\omega_K)^{-D}\prod_{j=1}^s W(\omega_{I^j})\in R(U(r))_{k,n}
\end{equation}
First using Equation ~\eqref{domingo}, we note $W(\omega_K)=1\times
x^{rn-r}$. Also $W(\omega_{I^j})=\bar{\lambda}(I^j)\times
x^{\codim(\omega_{I^j})}$. Then the $U(1)$ component of the product
in $\tilde{R}$ corresponds to $x$ raised to the number (mod $rn$)
$$\sum_j\codim(I^j)+Dr =r(n-r).$$

Therefore the coefficient of $W(\omega_{\{1,\dots,r\}})$ in the
product
$$W(\omega_K)^{-D}\prod_{j=1}^s W(\omega_{I^j})$$
equals the coefficient of the identity representation in the product
 $\prod_{j=1}^s\bar{\lambda}(I^j)$ in the Verlinde algebra of  $\SU(r)$ at level $\leq k$ ;  which equals
the right hand side of Equation ~\eqref{towit} as desired.
\end{proof}

\appendix\section{Some results on vector bundles on curves and Gromov-Witten theory}
\subsection{Schubert cycles}\label{fultone} Given  a subset $I$ of
$[n]$ of cardinality $r$ we will assume that it is written in the
form $\{i_1<i_2<\dots<i_r\}$. Let
$$E_{\sssize{\bullet}}:\text{ }\{0\}=E_0\subsetneq E_1\subsetneq\dots\subsetneq E_n=W$$
be a complete flag in the $n$-dimensional vector space $W$. Define
the Schubert cell $\Omega^o_{I}(E_{\sssize{\bullet}})\subseteq
\Gr(r,W)$ by
$$\Omega^o_I(E_{\sssize{\bullet}})=\{V\in \Gr(r,W)\mid \rk(V\cap E_u)=a \text{ for } i_a\leq u< i_{a+1},\text{ }a=0,\dots,r\}$$
where $i_0$ is defined to be $0$ and $i_{r+1}=n$.
$\Omega^o_{I}(E_{\sssize{\bullet}})$ is smooth. Its closure will  be
denoted by  ${\Omega}_{I}(E_{\sssize{\bullet}})$, and the cycle
class of this subvariety (in $H^*\Gr(r,W)$) is denoted by
$\omega_I$. For a fixed complete flag on $W$, it is easy to see that
every $r$-dimensional vector subspace belongs to a unique Schubert
cell.

The dual of $\omega_I$ under the intersection pairing is
$\omega_{I'}$ where $I'=\{n+1-i,i\in I\}$. This means that if
$\codim(\omega_I)+\codim(\omega_J)=r(n-r)$, the intersection number
$\omega_I\cdot\omega_J$ in $H^{2r(n-r)} X=\Bbb{Z}$ is $1$ if $J=I'$
and $0$ otherwise.
\subsection{Vector bundles on curves}
Let $X$ be a smooth curve of  genus $g$ and $E,F$ be vector bundles
of ranks $r$ and $k$ respectively such that $\chi(X,E\tensor F))=0$.
The following lemma is due to Faltings (see e.g. ~\cite{seshadri}).
\begin{lemma}\label{nori} If
$H^0(X,E\tensor F)=0$ then $E$ and $F$ are semi-stable.
\end{lemma}
\begin{lemma}\label{deform}
Let $X\to S$ be a projective, flat family of curves and $s_0\in S$.
Suppose $E_0$ is a vector bundle on $X_{s_0}$. Then, there exists an
\'{e}tale map $(T,t_0)\to (S,s_0)$ and  vector bundle $E_T$ on $X_T$
so that $E_{T,t_0}=E_{0}$. If $\det(E_0)$ is trivial then there
exist such $(T,t_0,E_T)$ so that $\det(E_T)$ is trivial.
\end{lemma}
\begin{proof} (Standard, for example see the proof of Proposition 4.1 on page 46 of ~\cite{SGA}.)
We  first  extend $E_0$ as a vector bundle. We will extend formally
first. For this (by Grothendieck's theorem on coherent sheaves on
 formal fibers), over successive thickenings of the central fibers.
Assume that $E$ been lifted to $\mathcal{X}_A\to \spec{A}$ (the base
change of $X\to S$ to $\spec{A}$) and $B$ is an Artinian local  ring
with an ideal $I$ such that $I^2=0$ and $B/I=A$. We want to extend
$E$ to $\mathcal{X}_B\to \spec(B)$. This is clear because (formally)
the obstruction to extension over Artin local rings is in
$H^2(\mathcal{X}_A, M_n({I}))$ ($n\times n$ matrices with entries in
${I}$) which vanishes on a curve. We  now apply Artin approximation
theorem ~\cite{artin} to obtain $(T,t_0)$.

If $\det(E_0)$ is trivial and we want the extension to preserve
this, the extension problem over successive Artin rings is
controlled by $H^2$ of traceless matrices with coefficients in
${I}$,  which is again $0$.
\end{proof}
\subsection{Shift operations in Gromov-Witten theory}\label{shiftie}
 Let $I^1,\dots,I^s$ be subsets of $[n]$ each of cardinality $r$ and $d$, $D$  integers. Suppose $I^{1}=\{i_1<\dots<i_r\}.$
Define $J$ a subset of $[n]$ of cardinality $r$ and an integer $\td$
as follows:
\begin{enumerate}
\item If $i_1>1$, let $J=\{i_1-1<\dots <i_r-1\}$ and $\tilde{d}=d$.
\item If $i_1=1$, let $J=\{i_2-1<\dots <i_r-1<n\}$ and $\tilde{d}=d-1$.
\end{enumerate}
The following proposition is proved in ~\cite{qh} (Proposition 2.5).
\begin{proposition}\label{shiftt1}
$$\langle \omega_{I^1},\omega_{I^2},\dots,\omega_{I^s}\rangle_{d,D}=\langle \omega_{J},\omega_{I^2},\dots,\omega_{I^s}\rangle_{\td,D-1}$$
\end{proposition}
We recall the reason for the equality (see ~\cite{qh} for more
details). We first verify that the expected dimensions of both the
intersections in the above equality of Gromov-Witten  numbers are
the same. Let $\mpp=\{{\pto}_1,\dots,{\pto}_s\}$ be a set of
(distinct points) on $\pone$ as before. Let $T$ be an ES bundle of
degree $-D$ and rank $n$. Choose generic complete flags
$E^{{\pto}_j}_{\bull}$ on the fibers $T_{{\pto}_j}$ for
$j=1,\dots,s$.

Let $\tilde{T}$ be the vector bundle which agrees with $T$ in
$\pone-\{{\pto}_1\}$ and whose sections in a small neighborhood $U$
of ${\pto}_1$ are sections $s$ of $T$ on $U-\{{\pto}_1\}$ such that
$ts$ is a holomorphic section of $T$ on $U$ whose fiber at
${\pto}_1$ lies in $E^{{\pto}_1}_1$ (the first element of the flag).
As coherent sheaves $\tilde{T}\supset T$. It can be shown that
$\tilde{T}$ is ES as well.

$\tilde{T}$ inherits complete flags from $T$ on its fibers at each
point of $\{{\pto}_1,\dots,{\pto}_s\}$. There is a bijection between
the set of subbundles of $T$ and those of $\tilde{T}$ both
restricted to $\pone-\{{\pto}_1\}$. Proposition ~\ref{shiftt1}
follows from this bijection and a calculation at ${\pto}_1$ (this
gives an inequality between the two Gromov-Witten numbers in
Proposition ~\ref{shiftt1} to start with, but this is a cyclic
process so we obtain equality.)

\bibliographystyle{plain}
\def\noopsort#1{}

\end{document}